\documentclass[preprint,sort&compress,final,10pt]{elsarticle}

\usepackage{amsmath}
\usepackage{amssymb}
\usepackage{graphicx}
\usepackage{latexsym}
\usepackage{epstopdf}
\usepackage{natbib}
\usepackage{color}
\usepackage{hyperref}

\newtheorem{theorem}{Theorem}[section]

\newtheorem{definition}[theorem]{Definition}
\newtheorem{lemma}[theorem]{Lemma}

\newtheorem{remark}[theorem]{Remark}
\numberwithin{equation}{section}

\def\Proof{\noindent{\bf Proof.}~}
\def\qed{\hfill$\square$\smallskip}

\def\dlim{\displaystyle\lim}

\def\Im{\mathrm{Im}}

\journal{\empty}
\date{}

\begin{document}

\begin{frontmatter}

\title{The reducibility of quasi-periodic linear Hamiltonian systems and its application to Hill's equation}

\author[au1,au2]{Nina Xue}

\address[au1]{School of Mathematical Sciences, Beijing Normal University, Beijing 100875, P.R. China.}

\author[au1]{Xiong Li\footnote{ Partially supported by the NSFC (11571041) and the Fundamental Research Funds for the Central Universities. Corresponding author.}}

\address[au2]{School of Mathematics and Information Sciences, Weifang University, Weifang, Shandong, 261061, P.R. China.}

\ead[au2]{xli@bnu.edu.cn}

\begin{abstract}
In this paper, we consider the reducibility of the quasi-periodic linear Hamiltonian system
$$\dot{x}=(A+\varepsilon Q(t))x, $$
where $A$ is a constant matrix with possible multiple eigenvalues, $Q(t)$ is analytic quasi-periodic with respect to $t$, and $\varepsilon$ is a sufficiently small parameter. Under some non-resonant conditions, it is proved that, for most sufficiently small $\varepsilon$, the Hamiltonian system can be reduced to a constant coefficient Hamiltonian system by means of a quasi-periodic symplectic change of variables with the same basic frequencies as $Q(t)$. Application to quasi-periodic Hill's equation is also given.
\end{abstract}

\begin{keyword}
Quasi-periodic Hamiltonian linear systems  \sep Reducibility \sep KAM iteration \sep  Hill's equation
\end{keyword}

\end{frontmatter}

\section{Introduction}

In this paper, we are concerned with the reducibility of the quasi-periodic linear Hamiltonian system
\begin{equation}\label{intro1}
\dot{x}=(A+\varepsilon Q(t))x,
\end{equation}
where $A$ is a constant matrix with possible multiple eigenvalues, $Q(t)$ is analytic quasi-periodic with respect to $t$, and $\varepsilon$ is a sufficiently small parameter.

Firstly, let us recall the definition of the reducibility for quasi-periodic linear systems. Let $A(t)$ be an $n\times n$ quasi-periodic matrix, the differential equation
\begin{equation}\label{GE}
\dot{x}=A(t)x, \ x\in \mathbb{R}^{n}
\end{equation}
is called reducible, if there exists a non-singular quasi-periodic change of variables
$$
x=\phi(t)y,
$$
where $\phi(t)$ and $\phi^{-1}(t)$ are quasi-periodic and bounded, which changes (\ref{GE}) into
\begin{equation}\label{GEC}
\dot{y}=By, \ y\in \mathbb{R}^{n}
\end{equation}
where $B$ is a constant matrix.

The well known Floquet theorem states that every periodic differential equation (\ref{GE}) can be reduced to a constant coefficient differential equation (\ref{GEC})  by means of a periodic change of variables with the same period as $A(t)$. However this is not true for the quasi-periodic linear system, one can see \cite{Palmer} for more details. In 1981, Johnson and Sell \cite{Johnson} proved that the  quasi-periodic linear system (\ref{GE}) is reducible if the quasi-periodic coefficient matrix $A(t)$ satisfies the "full spectrum" condition.

Therefore, many authors (\cite{Li1}, \cite{Jorba1}, \cite{Jorba2}, \cite{Xu1}) paid attention to the reducibility of the quasi-periodic linear system (\ref{intro1}),  which is close to a constant coefficient linear system. This problem was first studied by Jorba and Sim\'{o} in \cite{Jorba1}. Suppose that $A$ is a constant matrix with different eigenvalues, they proved that if the eigenvalues of $A$ and the frequencies of $Q$ satisfy some non-resonant conditions, then there exist some sufficiently small $\varepsilon_{0}>0$ and a non-empty Cantor set $E\subset(0,\varepsilon_{0})$, such that for any $|\varepsilon| \in E$, system \eqref{intro1} is reducible. Moreover, the relative measure of the set $(0,\varepsilon_{0})\setminus E$ in $(0,\varepsilon_{0})$ is exponentially small in $\varepsilon_{0}$. Later, Xu \cite{Xu1} obtained the similar result for the multiple eigenvalue case.

In 1996, Jorba and Sim\'{o} \cite{Jorba2} extended the conclusion of the linear system to the nonlinear system
 \begin{equation}\label{intro2}
\dot{x}=(A+\varepsilon Q(t,\varepsilon))x+\varepsilon g(t)+h(x,t), \ \ x\in\mathbb{R}^{n}.
\end{equation}
Suppose that $A$ has $n$ different nonzero eigenvalues, they proved that under some non-resonant conditions and non-degeneracy conditions, there exists a non-empty Cantor set $E\subset(0,\varepsilon_{0})$, such that for all $|\varepsilon| \in E$, system \eqref{intro2} is reducible. Later, Wang and Xu \cite{Wang} further inverstigated the nonlinear quasi-periodic system
 \begin{equation}\label{intro3}
\dot{x}=Ax+f(t,x,\varepsilon), \ \ x\in\mathbb{R}^{2},
\end{equation}
where $A$ is a real $2\times 2$ constant matrix, and $f(t,0,\varepsilon)=O(\varepsilon), \ \partial_x f(t,0,\varepsilon)=O(\varepsilon)$ as $\varepsilon\rightarrow 0.$
 They proved without any non-degeneracy condition, one of two results holds: (1) system \eqref{intro3} is reducible to $\dot{y}=By+O(y)$ for all $\varepsilon \in(0,\varepsilon_{0})$; (2) there exists a non-empty Cantor set $E\subset(0,\varepsilon_{0})$, such that system \eqref{intro3} is  reducible to $\dot{y}=By+O(y^{2})$ for all $|\varepsilon| \in E$.

 In \cite{You}, Her and You considered one-parameter family of quasi-periodic linear system
 \begin{equation}\label{You}
 \dot{x}=(A(\lambda)+g(\omega_1t, \cdots, \omega_lt, \lambda)x,
 \end{equation}
 where $A\in C^\omega(\Lambda, gl(m,C))$ ($C^\omega(\Lambda, gl(m,C))$ be the set of $m\times m$ matrices $A(\lambda)$ depending analytically on a parameter $\lambda$ in a closed interval $\Lambda\subset \mathbb{R}$), and $g$ is analytic and small. They proved that under some non-resonance conditions and non-degeneracy conditions, there exists an open and dense set $\mathcal{A}$ in $C^\omega(\Lambda, gl(m,C))$, such that for each $A\in \mathcal{A}$, system \eqref{You} is reducible for almost all $\lambda \in \Lambda$.

 Instead of a total reduction to a constant coefficient linear system, Jorba, Ramirez-ros and Villanueva  \cite{Jorba3} investigated the effective reducibility of the following quasi-periodic system
\begin{equation}\label{intro4}
\dot{x}=(A+\varepsilon Q(t,\varepsilon))x, \ \ |\varepsilon|\leq\varepsilon_{0},
\end{equation}
where $A$ is a constant matrix with different eigenvalues. They proved that under non-resonant conditions, by a quasi-periodic transformation,  system \eqref{intro4} is  reducible to a quasi-periodic system
$$\dot{y}=(A^{\ast}(\varepsilon)+\varepsilon R^{\ast}(t,\varepsilon))y, \ \ |\varepsilon|\leq \varepsilon_{\ast}\leq \varepsilon_{0},$$
where $R^{\ast}$ is exponentially small in $\varepsilon$. Li and Xu \cite{Li2} obtained the similar result for Hamiltonian systems.

In this paper we will study  the reducibility of the quasi-periodic linear {\it Hamiltonian} system (\ref{intro1}), where the matrix $A$ may have multiple eigenvalues. To this end, the following assumptions are made.

{\it Assumption A: Non-resonant condition.} Let all eigenvalues of the matrix $A$ be $\lambda_1, \cdots, \lambda_n$,  $Q(t)$ be an analytic quasi-periodic function on $D_{\rho}=\{\theta \in \mathbb{C}^{r}: |\Im\theta_{j}|\leq\rho, j=1,2,\cdots,r\}$ with the frequencies $\omega=(\omega_1, \cdots, \omega_r)$. Suppose that $\lambda=(\lambda_{1},  \cdots, \lambda_{n})$ and $\omega=(\omega_{1},  \cdots, \omega_{r})$ satisfy the non-resonant conditions
$$\big|\langle k,\omega\rangle\sqrt{-1}-\lambda_{i}+\lambda_{j}\big|\geq \frac{\alpha}{|k|^{\tau}}$$
for all $k\in \mathbb{Z}^{r}\setminus\{0\}$, $0\le i,j\le n$,  where $\alpha >0$ is a small constant and $\tau>r-1$.

{\it Assumption B: Non-degeneracy condition.} Assume that $A+\varepsilon \overline{Q}$ has $n$ different eigenvalues $\mu_{1},\cdots,\mu_{n}$ with
$|\mu_i|\geq 2\delta \varepsilon, \ |\mu_{i}-\mu_{j}|\geq 2\delta \varepsilon, \ i\not=j, \ 0\leq i,j\leq n,$ where $\delta$
is a positive constant independently of $\varepsilon$. Here we denote the average of $Q(t)$ by $\overline{Q}$, that is,
$$\overline{Q}=\lim_{T\rightarrow \infty}\frac{1}{2T}\int_{-T}^{T}Q(t)dt.$$

We are in a position to state the main result.

\begin{theorem}\label{Main result}

Suppose that the Hamiltonian system \eqref{intro1} satisfies the assumptions A and B. Then there exist some sufficiently small $\varepsilon_0>0$ and a non-empty Cantor subset $E_{\varepsilon_0}\subset (0,\varepsilon_0)$ with positive Lebesgue measure, such that for $\varepsilon \in E_{\varepsilon_0}$, the Hamiltonian system \eqref{intro1} is reducible, i.e.,
there is an analytic quasi-periodic \textbf{symplectic} transformation $x=\psi(t)y$, where $\psi(t)$ has same frequencies as $Q(t)$, which changes \eqref{intro1} into the Hamiltonian system $\dot{y}=By,$ where $B$ is a constant matrix. Moreover, if $\varepsilon_0$ is small enough, the relative measure of $E_{\varepsilon_0}$ in $(0,\varepsilon_0)$ is close to $1$.
\end{theorem}

Now we give some remarks on this result. Firstly, here we deal with the Hamiltonian system and have to find the \textit{symplectic} transformation, which is different from that in \cite{Jorba1} and \cite{Xu1}. Secondly, we consider the reducibility, other than the effective reducibility in \cite{Jorba3} and \cite{Li2}. The last but not the least, we can allow the matrix $A$ to have multiple eigenvalues. Of course, if the eigenvalues of $A$ are different, the non-degeneracy condition holds naturally.

After finishing this work, we consult references again and find the literature \cite{Li3}. In \cite{Li3}, Li, Zhu and Chen considered the following nonlinear analytic quasi-periodic Hamiltonian system
 \begin{equation}\label{Li}
 \dot{x}=(A+\varepsilon Q(t))x+\varepsilon g(t)+h(x,t), \ x\in \mathbb{R}^{2n},
 \end{equation}
where $A$ is a constant matrix, $h=O(x^2) (x\rightarrow 0)$, and $h(x,t), Q(t), g(t)$ are analytic quasi-periodic on $D_\rho$ with respect to t. They proved that, under suitable hypothesis of analyticity, non-resonance conditions and non-degeneracy conditions, there exists a non-empty Cantor set $E^\ast\subset(0,\varepsilon_{0})$ with positive Lebesgue measure, such that for  $\varepsilon\in E^\ast$, there is a quasi-periodic symplectic transformation, which changes the Hamiltonian system \eqref{Li} into the Hamiltonian system
$$\dot{y}=B(\varepsilon)y+h_\infty(y,t,\varepsilon),$$
where $B$ is a real constant matrix and $h_\infty(y,t,\varepsilon)=O(y^2)$ as $y\rightarrow 0$. Moreover, $\mbox{meas}((0,\varepsilon_{0})\setminus E^\ast)=o(\varepsilon_0)$ as $\varepsilon_0\rightarrow 0$.

Here we remark that if $g(t)\equiv 0, \ h(x,t)\equiv 0$ in \eqref{Li}, the result in Theorem 2.1 of \cite{Li3} is just the same as our main result in Theorem \ref{Main result}. However, in Theorem 2.1 of \cite{Li3}, $A$ is a matrix that can be diagonalized. In this paper, $A$ is only a constant matrix with possible multiple eigenvalues, which enables us to study equation \eqref{intro+1}, because
 \begin{equation*}
A=\left(
\begin{matrix}
0&1\\
0&0
\end{matrix}
\right)
\end{equation*}
can not be diagonalized. Moreover, the non-resonance conditions and non-degeneracy conditions in Theorem 2.1 of \cite{Li3} are all stronger than assumptions A and B in this paper. Of course, the proof of our main result in Theorem \ref{Main result} is different from that in \cite{Li3} in some respect. For instance, in the estimate on the measure, we do not need the non-degeneracy conditions which guarantee that $\lambda_i^m-\lambda_j^m$ are Lipschitz from above and from below. Furthermore, when proving the convergence of the iteration, our method can obtain some information to analyze the quasi-periodic Hill's equation \eqref{intro+1}. From the above, it is necessary to give the complete proof of Theorem \ref{Main result}. Therefore, we will prove Theorem \ref{Main result} in Section 3 of this paper.

As an example, we apply Theorem \ref{Main result} to the following quasi-periodic Hill's equation
\begin{equation}\label{intro+1}
\ddot{x}+\varepsilon a(t)x=0,
\end{equation}
where $a(t)$ is analytic quasi-periodic with the frequencies $\omega=(\omega_1, \cdots, \omega_r)$. Denote the average of $a(t)$ by $\bar{a}$.   If $\bar{a}>0$ and the frequencies $\omega$ of $a(t)$ satisfy the Diophantine condition
$$\big|\langle k,\omega\rangle\big|\geq \frac{\alpha}{|k|^{\tau}}$$
for all $k\in \mathbb{Z}^{r}\setminus\{0\}$, where $\alpha >0$ is a small constant and $\tau>r-1$, then there exists some sufficiently small $\varepsilon_0>0$,  equation \eqref{intro+1} is reducible and the equilibrium of \eqref{intro+1} is stable in the sense of Lyapunov for most sufficiently small $\varepsilon\in (0, \varepsilon_0)$. Moreover, all solutions of equation \eqref{intro+1} are quasi-periodic with the frequencies $\Omega=(\omega_1, \cdots,\omega_r, \sqrt{b})$ for most sufficiently small $\varepsilon\in (0, \varepsilon_0)$, where $b=\bar{a}\varepsilon+O(\varepsilon^2)$ as $\varepsilon\to 0$. Here we remark that if we rewrite equation \eqref{intro+1} into the Hamiltonian system \eqref{intro1}, we find that
\begin{equation*}
A=\left(
\begin{matrix}
0&1\\
0&0
\end{matrix}
\right),
\end{equation*}
which has multiple eigenvalues $\lambda_1=\lambda_2=0$. One can see Section 4 for more details about this example.

There are plenty of works about the stability of the equilibria of quasi-periodic Hamiltonian systems, one can refer to \cite{Bibikov}, \cite{Liu1}, \cite{Liu2} and \cite{Wu} for a detailed description. In general, in order to determine the type of stability of the equilibria of quasi-periodic Hamiltonian systems, the authors need to assume that the corresponding linearized system is reducible, and some conditions were added to the system after the reducibility. However, as far as we know, the case that the conditions are added to the original system has not been considered in the literature up to now, which we will study in the future.

The paper is organized as follows. In Section 2, we list some basic definitions and results that will be useful in the proof of the main result. In Section 3, we will prove Theorem \ref{Main result}. The quasi-periodic Hill's equation \eqref{intro+1} will be analyzed in Section 4.

\section{Some preliminaries}

We first give the definition of quasi-periodic functions.

\begin{definition}
A function $f$ is said to be a quasi-periodic function with a vector of basic frequencies $\omega=(\omega_{1},\omega_{2}, \cdots, \omega_{r})$, if $f(t)=F(\theta_{1},\theta_{2},\cdots,\theta_{r})$, where $F$ is $2\pi$ periodic in all its arguments and $\theta_{j}=\omega_{j}t$ for $j=1,2,\cdots,r.$ Moreover, if $F(\theta)\,(\theta=(\theta_{1},\theta_{2},\cdots,\theta_{r}))$ is analytic on $D_{\rho}=\{\theta \in \mathbb{C}^{r}:
|\Im\theta_{j}|\leq\rho, j=1,2,\cdots,r\}$, we say that $f(t)$ is analytic quasi-periodic on $D_{\rho}$.
\end{definition}

It is well known that an analytic quasi-periodic function $f(t)$ can be expanded as Fourier series
$$f(t)=\sum_{k\in\mathbb{Z}^{r}}f_{k}e^{\langle k,\omega\rangle\sqrt{-1}t}$$
with Fourier coefficients defined by
$$f_{k}=\frac{1}{(2\pi)^{r}}\int_{\mathbb{T}^{r}}F(\theta)e^{-\langle k,\theta\rangle\sqrt{-1}}d\theta.$$
Denote by $||f||_{\rho}$ the norm
$$||f||_{\rho}=\sum_{k\in\mathbb{Z}^{r}}|f_{k}|e^{|k|\rho}.$$

\begin{definition}
An $n\times n$ matrix $Q(t)=(q_{ij}(t))_{1\leq i,j\leq n}$ is said to be analytic quasi-periodic on $D_{\rho}$ with frequencies $\omega=(\omega_{1},\omega_{2}, \cdots, \omega_{r})$, if all $q_{ij}(t) \ (i,j=1,2,\cdots,n)$ are analytic quasi-periodic on $D_{\rho}$ with frequencies $\omega=(\omega_{1},\omega_{2}, \cdots, \omega_{r})$.
\end{definition}

Define the norm of $Q$ by
$$||Q||_{\rho}=\max_{1\leq i\leq n}\sum_{j=1}^{n}||q_{ij}||_{\rho}.$$
It is easy to see that
$$
||Q_{1}Q_{2}||_{\rho}\leq ||Q_{1}||_{\rho}||Q_{2}||_{\rho}.
$$
If $Q$ is a constant matrix, write $||Q||=||Q||_{\rho}$ for simplicity. Denote the average of $Q(t)$ by $\overline{Q}=(\overline{q}_{ij})_{1\leq i,j\leq n}$, where
$$\overline{q}_{ij}=\lim_{T\rightarrow \infty}\frac{1}{2T}\int_{-T}^{T}q_{ij}(t)dt,$$
see \cite{Siegel} for the existence of the limit.

Also we need two lemmas which are provided in this section for the proof of Theorem \ref{Main result}, that were proved in \cite{Jorba2}.

\begin{lemma}\label{lem1}
Let $h: B_\sigma(0)\subset \mathbb{R}^n\rightarrow \mathbb{R}^n$ be a $C^2$ function that satisfies $h(0)=0, D_xh(0)=0,$ $||D_{xx}h(x)||\leq K, x\in B_\sigma(0).$ Then $||h(x)||\leq\frac{K}{2}||x||^2$, $||D_xh(x)||\leq K||x||$.
\end{lemma}

\begin{lemma}\label{lem2}
Suppose that $B_0$ is an $n\times n$ matrix with different nonzero eigenvalues $\mu_1^0,\cdots,\mu^0_{n}$ satisfying $|\mu^0_i|>\gamma, \  |\mu_i^0-\mu_j^0|>\gamma, i \not=j, 1\leq i,j\leq n,$ and $S_0$ is a regular matrix such that $S_0^{-1}B_0S_0=diag(\mu_1^0,\cdots,\mu_{n}^0).$  Set  $\beta_0=\max\{||S_0||,||S_0^{-1}||\}$, and choose $b$ such that
$$0<b<\frac{\gamma}{(3n-1)\beta_0^2}.$$
If $B_1$ verifies $||B_1-B_0||\leq b,$ then the following conclusions hold:

(1) $B_1$ has n different nonzero eigenvalues $\mu_{1}^1,\cdots,\mu_{n}^1$.

(2) There exists a regular matrix $S_1$ such that $S_1^{-1}B_1S_1=diag(\mu_{1}^1,\cdots,\mu_{n}^1)$, which satisfies $||S_1||,
||S_1^{-1}||\leq \beta_1,$ where $\beta_1=2\beta_0$.
\end{lemma}

The next lemma is used to perform a step of the inductive procedure in the proof of Theorem \ref{Main result}.
\begin{lemma}\label{lem3}
Consider the differential equation
\begin{equation}\label{lemma1}
\dot{P}(t)=\Lambda P(t)-P(t)\Lambda+R(t),
\end{equation}
where $\Lambda$ is a constant Hamiltonian matrix with $n$ different eigenvalues $\nu_1,\cdots,\nu_n$, $R$ is an analytic quasi-periodic Hamiltonian matrix on $D_\rho$ with frequencies $\omega$, satisfying $\overline{R}=0$.

If
\begin{equation}\label{lemma2}
|\langle k,\omega\rangle\sqrt{-1}-\nu_i+\nu_j|\geq\frac{\alpha}{|k|^{3\tau}}
\end{equation}
for all $0\not=k\in\mathbb{Z}^r,$ and $|\nu_i|\geq \delta\varepsilon, \ |\nu_{i}-\nu_{j}|\geq \delta \varepsilon,$ for $i\not=j, \ 0\leq i,j\leq n,$ where $\delta$
is a positive constant independently of $\varepsilon$, then equation \eqref{lemma1} has a unique analytic quasi-periodic Hamiltonian solution $P(t)$ with $\overline{P}=0$, where $P(t)$ has frequencies $\omega$, and satisfies
\begin{equation}\label{lemma3}
||P||_{\rho-s}\leq\frac{c}{\alpha s^\nu }||R||_\rho
\end{equation}
with $\nu=3\tau+r$ and $0<s<\rho,$ where the constant $c$ depends only on $\tau$ and $r$.
\end{lemma}

\Proof Choosing $S$ such that $S^{-1}\Lambda S=D=diag(\nu_1,\cdots,\nu_n)$, making the change of variable $P(t)=SX(t)S^{-1}$ and defining $Y(t)=S^{-1}R(t)S$, equation \eqref{lemma1} becomes
 \begin{equation}\label{lemma+1}
\dot{X}(t)=DX(t)-X(t)D+Y(t), \ \ \overline{Y}=0.
\end{equation}
Expanding $X$ and $Y$ into Fourier series yields that
$$X(t)=\sum_{k\in \mathbb{Z}^r}X_ke^{\langle k,\omega\rangle\sqrt{-1}t}, \ \  Y(t)=\sum_{k\in \mathbb{Z}^r}Y_ke^{\langle k,\omega\rangle\sqrt{-1}t},$$
where $X_{k}=(x_{ij}^{k})_{1\leq i,j\leq n}$ and $Y_{k}=(y_{ij}^{k})_{1\leq i,j\leq n}.$

By comparing the coefficients of \eqref{lemma+1}, we obtain that
$$ x^k_{ij}=\frac{y^k_{ij}}{\langle k,\omega\rangle\sqrt{-1}-\nu_i+\nu_j}, \ \ 1\leq i,j\leq n, k\not=0, $$
and $$x^0_{ij}=0, \ 1\leq i,j\leq n.$$
Since $R$ is analytic on $D_\rho$, $Y$ is also analytic on $D_\rho$. Therefore, we have
$$||Y_k||\leq ||Y||_\rho e^{-|k|\rho}.$$
Hence
\begin{eqnarray*}\label{lemma+2}
||X||_{\rho-s}&=&\sum_{k\in\mathbb{Z}^r}||X_k||e^{|k|(\rho-s)}\\
&\leq&\sum_{0\not=k\in\mathbb{Z}^r}\frac{|k|^{3\tau} e^{-s|k|}}{\alpha}||Y||_\rho\\
&\leq&\frac{c}{\alpha s^\nu }||Y||_\rho,
\end{eqnarray*}
where $\nu=3\tau+r$ and $0<s<\rho.$ Here and hereafter we always use the same symbol $c$ to denote different constants in estimates.
Hence
$$||P||_{\rho-s}\leq c ||X||_{\rho-s}\leq\frac{c}{\alpha s^\nu }||Y||_\rho\leq \frac{c}{\alpha s^\nu }||R||_\rho.$$

Now we prove that $P$ is Hamiltonian. Since $\Lambda$ and $R$ are Hamiltonian, then $\Lambda=J\Lambda_J$ and $R=JR_J$,  where $\Lambda_J$ and $R_J$ are symmetric. Let $P_J=J^{-1}P$, if $P_J$ is symmetric, then $P$ is Hamiltonian. Below we prove that $P_J$ is symmetric. Substituting $P=JP_I$ into equation \eqref{lemma1} yields that
\begin{equation}\label{lemma4}
\dot{P_J}=\Lambda_JJP_J-P_JJ\Lambda_J+R_J,
\end{equation}
and transposing equation \eqref{lemma4}, we get
\begin{equation*}
\dot{P_J}^T=\Lambda_JJP_J^T-P_J^TJ\Lambda_J+R_J.
\end{equation*}
It is easy to see that $JP_J$ and $JP_J^T$ are solutions of \eqref{lemma1}, moreover, $\overline{JP_J}=\overline{JP_J^T}=0.$
Since the solution of \eqref{lemma1} with $\overline{P}=0$ is unique, we have that $JP_J=JP_J^T$, which implies that $P$ is Hamiltonian. Up to now, we have finished the proof of this lemma.\qed

\section{Proof of Theorem \ref{Main result}}

From the assumptions of Theorem \ref{Main result}, it follows that $A+\varepsilon \overline{Q}$ is a Hamiltonian matrix with $n$ different eigenvalues $\mu_{1}, \ \cdots, \ \mu_{n}$,
and $|\mu_i|\geq 2\delta \varepsilon, \ |\mu_{i}-\mu_{j}|\geq 2\delta \varepsilon, \ i\not=j, \ 0\leq i,j\leq n,$ where $\delta$
is a positive constant independently of $\varepsilon$.
We rewrite the Hamiltonian system \eqref{intro1} into
\begin{equation}\label{pf1}
\dot{x}=[A+\varepsilon \overline{Q}+\varepsilon(Q(t)-\overline{Q})]x:=(A_1+\varepsilon \widetilde{Q}(t))x,
\end{equation}
where $A_1=A+\varepsilon \overline{Q}, \ \widetilde{Q}(t)=Q(t)-\overline{Q}$, $\overline{\widetilde{Q}}=0$.

Introduce the change of variables $x=e^{\varepsilon P(t)}x_1$, where $P(t)$ will be determined later, under this symplectic transformation, the Hamiltonian system \eqref{pf1} is changed into the new Hamiltonian system
\begin{equation}\label{pf2}
\dot{x_1}=e^{-\varepsilon P(t)}(A_1+\varepsilon \widetilde{Q}-\varepsilon\dot{P})e^{\varepsilon P(t)}x_1.
\end{equation}
Expand $e^{\varepsilon P}$ and $e^{-\varepsilon P}$ into
$$e^{\varepsilon P}=I+\varepsilon P+B, \ \  e^{-\varepsilon P}=I-\varepsilon P+\widetilde{B},$$
where $$B=\frac{(\varepsilon P)^2}{2!}+\frac{(\varepsilon P)^3}{3!}+\cdots, \ \ \widetilde{B}=\frac{(\varepsilon P)^2}{2!}-\frac{(\varepsilon P)^3}{3!}+\cdots.$$
Then the Hamiltonian system \eqref{pf2} can be rewritten
\begin{eqnarray}\label{pf3}
\dot{x_1}&=&(I-\varepsilon P+\widetilde{B})(A_1+\varepsilon\widetilde{Q}-\varepsilon\dot{P})(I+\varepsilon P+B)x_1\nonumber\\[0.2cm]
&=&(A_1+\varepsilon \widetilde{Q}-\varepsilon \dot{P}+\varepsilon A_1P-\varepsilon PA_1+\varepsilon^2Q_1)x_1,
\end{eqnarray}
where
\begin{eqnarray*}
Q_1&=&-P(\widetilde{Q}-\dot{P})+(\widetilde{Q}-\dot{P})P-P(A_1+\varepsilon\widetilde{Q}-\varepsilon\dot{P})P\\
&&+(I-\varepsilon P)(A_1+\varepsilon\widetilde{Q}-\varepsilon\dot{P})\frac{B}{\varepsilon^2}\\
&&+e^{\varepsilon P}(A_1+\varepsilon\widetilde{Q}-\varepsilon\dot{P})\frac{\widetilde{B}}{\varepsilon^2}.
\end{eqnarray*}

We would like to have
$$ \widetilde{Q}-\dot{P}+ A_1P-PA_1=0,$$
which is equivalent to
\begin{equation}\label{pf4}
\dot{P}=A_1P-PA_1+\widetilde{Q}.
\end{equation}

By the assumption B of Theorem \ref{Main result}, it is easy to see that the inequalities
$$|\mu_i|\geq \delta\varepsilon, \ |\mu_{i}-\mu_{j}|\geq \delta \varepsilon, \ i\not=j, \ 0\leq i,j\leq n$$
hold. Moreover, if the equalities
\begin{equation}\label{+1}
|\langle k,\omega\rangle\sqrt{-1}-\mu_i+\mu_j|\geq\frac{\alpha_0}{|k|^{3\tau}}, \ 0\not=k\in\mathbb{Z}^r,
\end{equation}
also hold, where $\alpha_0=\frac{\alpha}{2},$ thus, by Lemma \ref{lem3}, \eqref{pf4} is solvable for $P$ on a smaller domain, that is, there is a unique quasi-periodic Hamiltonian matrix $P(t)$ with frequencies $\omega$ on $D_{\rho-s}$, which satisfies $\overline{P}=0$ and
\begin{equation}\label{pf5}
||P||_{\rho-s}\leq\frac{c}{\alpha_0 s^\nu}||\widetilde{Q}||_{\rho}\leq\frac{c}{\alpha_0 s^\nu}||Q||_{\rho},
\end{equation}
where $ s=\frac{1}{2}\rho.$

Therefore, by \eqref{pf4},  the Hamiltonian system \eqref{pf3} becomes
\begin{equation}\label{pf5}
\dot{x_1}=(A_1+\varepsilon^2Q_1)x_1,
\end{equation}
where
\begin{eqnarray*}
Q_1&=&P(A_1P-PA_1)+(PA_1-A_1P)P\\[0.2cm]
&&-P(A_1+\varepsilon(PA_1-A_1P))P\\
&&+(I-\varepsilon P)(A_1+\varepsilon(PA_1-A_1P))\frac{B}{\varepsilon^2}\\
&&+\frac{\widetilde{B}}{\varepsilon^2}(A_1+\varepsilon(PA_1-A_1P))e^{\varepsilon P}.
\end{eqnarray*}

From Lemma \ref{lem1}, it follows that
$$||B||_{\rho-s}\leq c||\varepsilon P||_{\rho-s}^2, \ ||\widetilde{B}||_{\rho-s}\leq c||\varepsilon P||_{\rho-s}^2.$$
Therefore, if $|\varepsilon|$ is sufficiently small, we have
$$||Q_1||_{\rho-s}\leq c||P||_{\rho-s}^2\leq\frac{c}{\alpha_0^2s^{2\nu}}||Q||_{\rho}^2.$$

Now we consider the iteration step. In the $m^{th}$ step, we consider the Hamiltonian system
\begin{equation}\label{pf6}
\dot{x}_{m}=(A_{m}+\varepsilon^{2^{m}} Q_{m}(t))x_{m}, \ m\geq 1,
\end{equation}
where $A_m$ has $n$ different eigenvalues $\lambda_1^m,\cdots,\lambda_n^m$ with
$$ |\lambda_i^m|\geq \delta \varepsilon,  \ \ |\lambda_i^m-\lambda_j^m|\geq \delta \varepsilon, \ i\not=j, \ 1\leq i,j\leq n,$$
here we define $\lambda_i^1=\mu_i, \ i=1,\cdots,n. $

Let $A_{m+1}=A_m+\varepsilon^{2^{m}}\overline{Q}_m$, then the Hamiltonian system \eqref{pf6} becomes
\begin{equation}\label{pf7}
\dot{x}_{m}=(A_{m+1}+\varepsilon^{2^{m}} \widetilde{Q}_{m})x_{m}, \ m\geq 1,
\end{equation}
where $\widetilde{Q}_m=Q_m(t)-\overline{Q}_m.$

We need to solve
$$\dot{P}_m=A_{m+1}P_m-P_mA_{m+1}+\widetilde{Q}_m.$$
If$$|\langle k,\omega\rangle\sqrt{-1}-\lambda_i^{m+1}+\lambda_j^{m+1}|\geq\frac{\alpha_{m}}{|k|^{3\tau}}, \ \  0\not=k\in\mathbb{Z}^r,$$
and
$A_{m+1}$ has $n$ different eigenvalues $\lambda_1^{m+1},\cdots,\lambda_n^{m+1}$ with
$$ |\lambda_i^{m+1}|\geq \delta \varepsilon, \ \
|\lambda_i^{m+1}-\lambda_j^{m+1}|\geq \delta \varepsilon,\ \ i\not=j,\ 1\leq i,j\leq n,$$
by Lemma \ref{lem3}, there is a unique quasi-periodic Hamiltonian matrix $P_m(t)$ with frequencies $\omega$ on $D_{\rho_m-s_m}$, which satisfies
\begin{equation}\label{pf8}
||P_m||_{\rho_m-s_m}\leq\frac{c}{\alpha_m s_m^\nu }||Q_m||_{\rho_m}.
\end{equation}
Thus, under the symplectic change of variables $x_{m}=e^{\varepsilon^{2^{m}}P_m(t)}x_{m+1}$, the Hamiltonian system \eqref{pf7} is changed into
\begin{equation*}
\dot{x}_{m+1}=(A_{m+1}+\varepsilon^{2^{m+1}}Q_{m+1})x_{m+1},
\end{equation*}
where
\begin{eqnarray*}
Q_{m+1}(t)&=&P_m(A_{m+1}P_m-P_m A_{m+1})+(P_m A_{m+1}-A_{m+1}P_m)P_m\\[0.2cm]
&&-P_m\left(A_{m+1}+\varepsilon^{2^{m}}(P_m A_{m+1}-A_{m+1}P_m)\right)P_m\\
&&+(I-\varepsilon^{2^{m}}P_m)\left(A_{m+1}+\varepsilon^{2^{m}}(P_m A_{m+1}-A_{m+1}P_m)\right)\frac{B_m}{\varepsilon^{2^{m+1}}}\\
&&+\frac{\widetilde{B}_m}{\varepsilon^{2^{m+1}}}\left(A_{m+1}+\varepsilon^{2^{m}}(P_m A_{m+1}-A_{m+1}P_m)\right)e^{\varepsilon^{2^{m}} P_m},
\end{eqnarray*}
$$e^{\varepsilon^{2^{m}} P_m}=I+\varepsilon^{2^{m}} P_m+B_m,$$ $$ \ e^{-\varepsilon^{2^{m}} P_m}=I-\varepsilon^{2^{m}} P_m+\widetilde{B}_m,$$
and
$$B_m=\frac{(\varepsilon^{2^{m}} P_m)^2}{2!}+\frac{(\varepsilon^{2^{m}}P_m)^3}{3!}+\cdots,$$
$$\widetilde{B}_m=\frac{(\varepsilon^{2^{m}} P_m)^2}{2!}-\frac{(\varepsilon^{2^{m}}P_m)^3}{3!}+\cdots.$$

 From Lemma \ref{lem1}, it follows that
$$||B_m||_{\rho_m-s_m}\leq c||\varepsilon^{2^m} P_m||_{\rho_m-s_m}^2, \ ||\widetilde{B}_m||_{\rho_m-s_m}\leq c||\varepsilon^{2^m} P_m||_{\rho_m-s_m}^2.$$
Therefore, if $|\varepsilon|$ is sufficiently small,
by \eqref{pf8} we have
\begin{equation}\label{pf9}
||Q_{m+1}||_{\rho_m-s_m}\leq\frac{c}{\alpha_m^2 s_m^{2\nu}}||Q_m||^2_{\rho_m}.
\end{equation}

Now we prove that the iteration is convergent as $m\rightarrow \infty$. When $m=1$, we choose
$$ \alpha_1=\frac{1}{4}\alpha, \  \rho_1=\frac{1}{2}\rho, \ s_1=\frac{1}{8}\rho, \ F_1=\frac{||\varepsilon^2 Q_1||_{\rho_1}}{\alpha_1^2s_1^{2\nu}}.$$
At the $m^{th}$ step, we define
$$\alpha_m=\frac{\alpha}{(m+1)^2}, \  s_m=\frac{\rho}{2^{m+2}}, \  \rho_m=\rho_1-(s_1+s_2+\cdots+s_{m-1})$$ and
$$F_m=\frac{\varepsilon^{2^{m}}||Q_m||_{\rho_m}}{\alpha_m^2s_m^{2\nu}}.$$
By \eqref{pf9}, we have
$$F_{m+1}\leq c\frac{\varepsilon^{2^{m+1}}||Q_m||^2_{\rho_m}}{(\alpha_m^2s_m^{2\nu})^2}=cF_m^2,$$
where the constant $c$ depends only on $\alpha, \rho$. Hence it follows that
\begin{equation}\label{pf10}
cF_{m+1}\leq(cF_m)^2\leq(cF_1)^{2^m}.
\end{equation}

If $cF_1<1$, then $cF_m\rightarrow 0$ as $m\rightarrow\infty.$ From \eqref{pf8} it follows that
\begin{equation}\label{pf11}
||\varepsilon^{2^{m}}P_m||_{\rho_m-s_m}<c F_m.
\end{equation}
Thus, if $cF_1<\frac{1}{2}$, then
$$||e^{\pm \varepsilon^{2^{m}}P_m}||_{\rho_m}\leq 2.$$
Since
\begin{equation}\label{pf+1}
||A_{m+1}-A_m||=||\varepsilon^{2^{m}}\overline{Q}_m||\leq ||\varepsilon^{2^{m}}Q_m||_{\rho_m}< c F_m,
\end{equation}
if $cF_1\leq\frac{\delta \varepsilon}{(3n-1)\beta_m^2}$, it follows from \eqref{pf+1} that
$$||A_{m+1}-A_m||\leq\frac{\delta \varepsilon}{(3n-1)\beta_m^2}, \ \mbox{for} \  \mbox{any} \ m\geq 1,$$
where $\beta_m=\max\{||S_m||,||S_m^{-1}||\}$ and $S_m$
is the regular matrix in Lemma \ref{lem2} such that
$$S_m^{-1}A_mS_m=diag(\lambda_1^m, \cdots, \lambda_n^m).$$
Thus, it follows from Lemma \ref{lem2} that
$A_{m+1}$ has $n$ different eigenvalues $\lambda_1^{m+1},\cdots,\lambda_n^{m+1}.$  Moreover
$$|\lambda_{i}^{m+1}-\lambda_{j}^{m+1}|\geq \delta \varepsilon, \ i\not=j, 1\leq i,j\leq n,$$
and
$$|\lambda_i^{m+1}|\geq \delta\varepsilon, \ i=1, \cdots, n.$$
In fact,
 \begin{eqnarray*}
\big|\lambda_{i}^{m+1}-\lambda_{j}^{m+1}\big|&\geq& \big|\lambda^1_{i}-\lambda^1_{j}\big|
-\sum_{l=1}^{m}\left(|\lambda_{i}^{l+1}-\lambda^l_{i}|+|\lambda_{j}^{l+1}-\lambda^l_{j}|\right)\\
&\geq&\big|\lambda^1_{i}-\lambda^1_{j}\big|-2\sum_{l=1}^{m}||A_{l+1}-A_l||\\
&\geq&\big|\lambda^1_{i}-\lambda^1_{j}\big|-2\sum_{l=1}^{m}c F_l\\
&\geq&2\delta \varepsilon-2\sum_{l=1}^{m}c F_l.
\end{eqnarray*}
Moreover, we have
\begin{eqnarray*}
\sum_{l=1}^{m}c F_l&\leq& \sum_{l=1}^{\infty}c F_l\leq \sum_{m=0}^{\infty}(c F_1)^{2^{m}}\leq\sum_{m=1}^{\infty}(c F_1)^m \\
&=&\frac{c F_1}{1-c F_1}<2c F_1.
\end{eqnarray*}
Thus, if $cF_1\leq \min\left\{\frac{1}{2},\frac{1}{4}\delta \varepsilon, \frac{\delta \varepsilon}{(3n-1)\beta_m^2}\right\}$, that is, $0< \varepsilon \leq \min\left\{1,\frac{c}{||Q||_\rho},\frac{c}{||Q||^2_\rho}\right\},$
then by \eqref{pf10}, we have
$$|\lambda_{i}^{m+1}-\lambda_{j}^{m+1}|\geq2\delta\varepsilon-4\varepsilon cF_1\geq \delta \varepsilon, \ i\not=j, 1\leq i,j\leq n.$$
In the same way as above, we have
$$|\lambda_i^{m+1}|\geq \delta\varepsilon, \ i=1, \cdots, n.$$

Let $D_\ast=\bigcap_{m=1}^\infty D_{\rho_m}=D_{\frac{\rho}{4}}.$ By \eqref{pf11}, the composition of all the changes $e^{\varepsilon^{2^{m}} P_m}$ converges to $\psi$ as $m\rightarrow\infty$.
Obviously, $$||\varepsilon^{2^{m}}Q_m||_{D_\ast}\leq cF_m\rightarrow 0, \ m\rightarrow\infty.$$
Furthermore, it follows from \eqref{pf+1} that $A_m$ is convergent as $m\rightarrow\infty$. Define $B=\dlim_{m\rightarrow\infty}A_m.$
Then, under the symplectic change of variables $x=\psi(t) y$, the Hamiltonian system \eqref{intro1} is changed into
$\dot{y}=B y.$

Now we prove that, for most sufficiently small $\varepsilon$, such symplectic transformation exists. From the above iteration, we need to prove that the non-resonant conditions
\begin{equation}\label{pf12}
|\langle k,\omega\rangle\sqrt{-1}-\lambda_i^{m+1}+\lambda_j^{m+1}|\geq\frac{\alpha_m}{|k|^{3\tau}}
\end{equation}
for all $0\not=k\in\mathbb{Z}^r, 1\leq i,j\leq n, \  m=0,1,2,\cdots$,  hold for most sufficiently small $\varepsilon$.

Let $f(\varepsilon)=\langle k,\omega\rangle\sqrt{-1}-\lambda_i^{m+1}+\lambda_j^{m+1}, \ i\not=j,$ and
$$O_{ijm}^k=\left\{\varepsilon\in(0,\varepsilon_0):|f(\varepsilon)|<\frac{\alpha_m}{|k|^{3\tau}}\right\},$$
where we choose
$$\varepsilon_0=\min\left\{1, \frac{c}{||Q||_\rho},\frac{c}{||Q||^2_\rho}\right\} $$  such that, for $\varepsilon\in(0,\varepsilon_0)$, the above iteration is convergent, and
\begin{equation}\label{pf13}
\Big|\frac{df}{d\varepsilon}\Big|=\Big|\frac{d}{d\varepsilon}(\lambda_i^{m+1}-\lambda_j^{m+1})\Big|\geq \delta.
\end{equation}
 For $\varepsilon \in (0,\varepsilon_0)$, by \eqref{pf+1} we have
\begin{eqnarray*}
||A_{m+1}-A_1||&\leq&||A_{m+1}-A_{m}||+\cdots+||A_2-A_1||\\[0.2cm]
&\leq&c F_m+\cdots+ c F_2 \\
&\leq& 2c F_1\leq \frac{1}{2}\delta\varepsilon.
\end{eqnarray*}
Hence
\begin{eqnarray*}
|f(\varepsilon)|&\geq&|\langle k,\omega\rangle\sqrt{-1}-\lambda_i+\lambda_j|-|\lambda_i^1-\lambda_i|-|\lambda_j^1-\lambda_j|\\[0.2cm]
&&-|\lambda_i^1-\lambda_i^{m+1}|-|\lambda_j^1-\lambda_j^{m+1}|\\
&\geq&\frac{\alpha}{|k|^\tau}-2q\varepsilon-2||A_{m+1}-A_1||\\
&\geq&\frac{\alpha}{|k|^\tau}-2q\varepsilon-\delta\varepsilon
\geq\frac{\alpha}{|k|^\tau}-3M\varepsilon_0,
\end{eqnarray*}
where $M=\max\{q,\delta\}$.

If $\frac{1}{|k|^\tau}\geq\frac{6 M\varepsilon_0}{\alpha},$ then
$$|f(\varepsilon)|\geq\frac{\alpha}{2|k|^\tau}>\frac{\alpha_m}{|k|^{3\tau}},$$
and $O_{ijm}^k=\emptyset.$

Suppose that $\frac{1}{|k|^\tau}<\frac{6 M\varepsilon_0}{\alpha}$. By \eqref{pf13}, it follows that
$$ \mbox{meas}\left(O_{ijm}^k\right)<\frac{\alpha_m}{|k|^{3\tau}\delta}.$$
Thus,
\begin{eqnarray*}
\mbox{meas}\left(\bigcup_{i\not=j}\bigcup_{0\not=k\in\mathbb{Z}^r}O_{ijm}^k\right)
&\leq&\frac{n^2\alpha_m}{\delta}\sum_{|k|^{-\tau}<\frac{6 M\varepsilon_0}{\alpha}}\frac{1}{|k|^{3\tau}}\\
&\leq& \frac{n^2\alpha_m}{\delta}\cdot\frac{36M^2\varepsilon_0^2}{\alpha^2}\sum_{k\in\mathbb{Z}^r}\frac{1}{|k|^\tau}\\
&\leq&\frac{c\varepsilon_0^2}{m^2}.
\end{eqnarray*}
Let
$$E_m=\left\{\varepsilon\in(0,\varepsilon_0):|\langle k,\omega\rangle\sqrt{-1}-\lambda_i^{m+1}+\lambda_j^{m+1}|>\frac{\alpha_m}{|k|^{3\tau}}, 0\not=k\in\mathbb{Z}^r, \ i\not=j\right\}.$$
Then $$(0,\varepsilon_0)-E_m=\bigcup_{i\not=j}\bigcup_{0\not=k\in\mathbb{Z}^r}O_{ijm}^k.$$
Thus
$$\mbox{meas}\left((0,\varepsilon_0)-E_m\right)\leq \frac{c\varepsilon_0^2}{m^2}.$$
Let $E_{\varepsilon_0}=\cap_{m=1}^{\infty}E_m,$ then $$\mbox{meas}\left((0,\varepsilon_0)-E_{\varepsilon_0}\right)\leq c\varepsilon_0^2,$$
and
$$\lim_{\varepsilon_0\rightarrow 0}\frac{\mbox{meas}\left((0,\varepsilon_0)-E_{\varepsilon_0}\right)}{\varepsilon_0}=0.$$
Therefore, $E_{\varepsilon_0}$ is a non-empty subset of $(0,\varepsilon_0)$. Thus, for $\varepsilon\in E_{\varepsilon_0}$, the Hamiltonian system \eqref{intro1} is reducible.
i.e., there exists a symplectic transformation $x=\psi(t)y$, which changes the Hamiltonian system \eqref{intro1} into the Hamiltonian system $\dot{y}=By$. Thus, Theorem \ref{Main result} is proved completely.

\section{The quasi-periodic Hill's equation}

 As an example, we apply Theorem \ref{Main result} to the following quasi-periodic Hill's equation
\begin{equation}\label{exp1}
\ddot{x}+\varepsilon a(t)x=0,
\end{equation}
where $a(t)$ is an analytic quasi-periodic function on $D_\rho$ with frequencies $\omega=(\omega_1, \cdots, \omega_r)$. Denote the average of $a(t)$ by $\bar{a}$, and suppose $\bar{a}>0$.

Let $\dot{x}=y$, then equation \eqref{exp1} can be rewritten in the equivalent form
\begin{equation}\label{exp2}
\dot{x}=y, \ \ \dot{y}=-\varepsilon a(t)x.
\end{equation}
To apply Theorem \ref{Main result}, we express \eqref{exp2} in the form
\begin{equation}\label{exp3}
\dot{z}=(A+\varepsilon Q(t))z,
\end{equation}
where
\begin{equation*}
z=\left(
\begin{matrix}
x\\
y
\end{matrix}
\right), \ \
A=\left(
\begin{matrix}
0&1\\
0&0
\end{matrix}
\right), \ \
Q=\left(
\begin{matrix}
0&0\\
-a(t)&0
\end{matrix}
\right).
\end{equation*}
 It is easy to see that $A$ has multiple eigenvalues $\lambda_1=\lambda_2=0$, moreover, $A+\varepsilon \overline{Q}$ has two different eigenvalues $\mu_1=i\sqrt{\bar{a}\varepsilon}, \ \mu_2=-i\sqrt{\bar{a}\varepsilon},$
where $ \overline{Q}$ stands for the average of the matrix $Q(t)$ and $i=\sqrt{-1}.$  It is clear that
$$|\mu_i|=\sqrt{\bar{a}}\sqrt{\varepsilon}\geq 2\delta \varepsilon, \ i=1,2,$$
and
$$|\mu_1-\mu_2|=2\sqrt{\bar{a}}\sqrt{\varepsilon}\geq 2\delta \varepsilon,$$
where we choose $\delta=\frac{1}{2}\sqrt{\bar{a}}>0$, which is a
constant independent of $\varepsilon$. Therefore, Theorem \ref{Main result} can be applied. It follows from Theorem \ref{Main result} that the following result holds.

\begin{theorem}\label{app1}
Assume that $a(t)$ is an analytic quasi-periodic function on $D_\rho$ with frequencies $\omega=(\omega_1, \cdots, \omega_r)$, and $\bar{a}>0$.
If the frequencies $\omega$ of $a(t)$ satisfy the Diophantine condition
\begin{equation}\label{Diophantine condition}
\big|\langle k,\omega\rangle\big|\geq \frac{\alpha}{|k|^{\tau}}
\end{equation}
for all $k\in \mathbb{Z}^{r}\setminus\{0\}$, where $\alpha >0$ is a small constant and $\tau>r-1$.

Then there exist some sufficiently small $\varepsilon_0>0$ and a non-empty Cantor subset $E_{\varepsilon_0}\subset (0,\varepsilon_0)$ with positive Lebesgue measure, such that for $\varepsilon \in E_{\varepsilon_0}$, the Hamiltonian system \eqref{exp3} is reducible. Moreover, if $\varepsilon_0$ is small enough, the relative measure of $E_{\varepsilon_0}$ in $(0,\varepsilon_0)$ is close to $1$.
\end{theorem}
\begin{remark}
From Theorem \ref{app1}, it follows that equation \eqref{exp1} can be changed into a constant coefficient differential equation for most sufficiently small $\varepsilon>0$.
\end{remark}

Now we want to study the Lyapunov stability of the equilibrium of the equation \eqref{exp1}, using the results obtained in Section 3.
If $a(t)$ is periodic in time ($T$ is the period),  one famous stability criterion was given by Magnus and Winkler \cite{Magnus} for Hill's equation
\begin{equation}\label{sta1}
\ddot{x}+ a(t)x=0.
\end{equation}
That is, \eqref{sta1} is stable if
$$a(t)>0, \ \int_0^Ta(t)dt\leq\frac{4}{T},$$
which can be shown using a Poincar\'{e} inequality. Such a stability criterion had been generalized and improved by Zhang and Li in \cite{Zhang1}, which now is the so-called $L^p$-criterion. Recently, Zhang in \cite{Zhang2} had extended such a criterion to the linear planar Hamiltonian system
$$\dot{x}=m(t)y, \ \ \dot{y}=-n(t)x,$$
where $m(t), n(t)$ are continuous and $T$-periodic functions.

However, for the quasi-periodic Hill's equation \eqref{exp1}, the results above can not be applied directly.
Now we obtain a result about the stability of the equilibrium of equation \eqref{exp1}.
\begin{theorem}\label{stable}
Under the conditions of Theorem \ref{app1}, the equilibrium of the equation \eqref{exp1}
is stable in the sense of Lyapunov for most sufficiently small $\varepsilon>0$.
\end{theorem}

\Proof Theorem \ref{app1} tells us that, for most sufficiently small $\varepsilon \in(0,\varepsilon_0)$, there exists an analytic quasi-periodic symplectic transformation
$z = \psi(t)z_\infty$, where $\psi(t)$ has same frequencies as $Q(t)$, which changes \eqref{exp3}
into the Hamiltonian system
\begin{equation}\label{stable+1}
\dot{z}_\infty= Bz_\infty,
\end{equation}
where $B$ is a constant matrix. Moreover, from the proof of Theorem \ref{Main result} in Section 3, it follows that
$B$ has two different eigenvalues $\lambda_1^\infty, \lambda_2^\infty$, satisfying
$$|\lambda_i^\infty|\geq \delta\varepsilon \ (i=1,2), \ \
|\lambda_1^\infty-\lambda_2^\infty|\geq \delta\varepsilon.$$
 Furthermore, by the proof of Theorem \ref{Main result}, we have
\begin{equation}\label{stable2}
||B-(A+\varepsilon \overline{Q})||\leq cF_1=O(\varepsilon^2).
\end{equation}
Therefore, the two different eigenvalues of $B$ are pure imaginary and can be written in the form
$$\lambda_i^\infty=\pm i\sqrt{b}, \ i=1,2,$$
where
$b$ can be written in the following form
\begin{equation}\label{b}
b=\bar{a}\varepsilon+O(\varepsilon^2),
\end{equation}
which depends on $\bar{a}$ and $\varepsilon$ only.

Thus, there exists a singular symplectic matrix $S$ such that
\begin{equation*}
S^{-1}BS=\left(
\begin{matrix}
i\sqrt{b}&0\\
0&-i\sqrt{b}
\end{matrix}
\right).
\end{equation*}
Let $z_\infty=S\tilde{z}_\infty$,  under this symplectic transformation, the Hamiltonian system \eqref{stable+1} is changed into
\begin{equation*}
\dot{\tilde{z}}_\infty=S^{-1}BS\tilde{z}_\infty=\left(
\begin{matrix}
i\sqrt{b}&0\\
0&-i\sqrt{b}
\end{matrix}
\right)\tilde{z}_\infty.
\end{equation*}
Hence, by an analytic quasi-periodic symplectic transformation, equation \eqref{exp1} is changed into
\begin{equation}\label{exp4}
\ddot{x}_\infty+b x_\infty=0.
\end{equation}
It is easy to see that equation \eqref{exp4} is elliptic. Therefore, the equilibrium of equation \eqref{exp1}
is stable in the sense of Lyapunov for most sufficiently small $\varepsilon>0$.\qed

For the existence of quasi-periodic solution of equation \eqref{exp1}, we have the following result.
\begin{theorem}\label{quasi-periodic solution}
Under the conditions of Theorem \ref{app1}, all solutions of equation \eqref{exp1} are quasi-periodic  with frequencies $\Omega=(\omega_1, \cdots,\omega_r, \sqrt{b})$ for most sufficiently small $\varepsilon>0$, where $b$ is given by (\ref{b}).
\end{theorem}
\Proof By Theorem \ref{app1},  we know that, for most sufficiently small $\varepsilon \in(0,\varepsilon_0)$, there exists an analytic quasi-periodic symplectic transformation
which has same frequencies as $a(t)$, by this transformation, equation \eqref{exp1} is changed into \eqref{exp4}. On the other hand,
it is easy to see that all solutions of the equation \eqref{exp4} are periodic, and the frequency of these solutions is $\sqrt{b}$.

Thus, we only need to prove that, for most sufficiently small $\varepsilon\in(0,\varepsilon_0)$, the following non-resonant condition
\begin{equation}\label{solution1}
|k_1\omega_1+\cdots+k_r\omega_r+k_{r+1}\sqrt{b}|\geq\frac{\alpha_0}{|k|^{5\tau+4}}
\end{equation}
holds for all $k=(k_1,\cdots,k_{r+1})\in\mathbb{Z}^{r+1}\setminus\{0\}$ , where $\alpha_0$ is defined in Section 3, that is, $\alpha_0=\frac{1}{2}\alpha$, and $\omega=(\omega_1,\cdots,\omega_r)$ are the frequencies of $a(t)$.

If $k_{r+1}=0$, then from the Diophantine condition \eqref{Diophantine condition},  it follows that \eqref{solution1} holds.

Suppose that $k_{r+1}\not=0$. Let $g(\varepsilon)=k_1\omega_1+\cdots+k_r\omega_r+k_{r+1}\sqrt{b},
$ and
$$O_k=\left\{\varepsilon\in(0,\varepsilon_0):|g(\varepsilon)|<\frac{\alpha_0}{|k|^{5\tau+4}}\right\}.$$
It follows from the non-degeneracy condition that
\begin{equation}\label{solution2}
\Big|\frac{dg}{d\varepsilon}\Big|=\Big|\frac{d}{d\varepsilon}(k_{r+1}\sqrt{b})\Big|\geq |k_{r+1}|\delta.
\end{equation}

By \eqref{b}, we have
$$\sqrt{b}\leq 4\delta\sqrt{\varepsilon}.$$
From the Diophantine condition \eqref{Diophantine condition}, it follows that
\begin{eqnarray*}
|g(\varepsilon)|&\geq&\frac{\alpha}{(|k_1|+\cdots+|k_r|)^\tau}-|k_{r+1}|\sqrt{b}\\
&\geq&\frac{\alpha}{|k|^\tau}-|k_{r+1}|\sqrt{b}\\
&\geq& \frac{\alpha}{|k|^\tau}-|k_{r+1}|4\delta\sqrt{\varepsilon}\\
&\geq&\frac{\alpha}{|k|^\tau}-4\delta|k|\sqrt{\varepsilon_0}.
\end{eqnarray*}

If $\frac{1}{|k|^{\tau+1}}\geq \frac{8\delta\sqrt{\varepsilon_0}}{\alpha}$, then
$$\big|g(\varepsilon)\big|\geq\frac{\alpha}{2|k|^\tau}\geq\frac{\alpha_0}{|k|^{5\tau+4}},$$
and $O_k=\emptyset$.

Suppose that $\frac{1}{|k|^{\tau+1}}< \frac{8\delta\sqrt{\varepsilon_0}}{\alpha}$, it follows from \eqref{solution2} that
$$\mbox{meas}(O_k)<\frac{\alpha_0}{|k|^{5\tau+4}|k_{r+1}|\delta}.$$
Thus,
\begin{eqnarray*}
\mbox{meas}\left(\bigcup_{0\not=k\in\mathbb{Z}^{r+1}}O_k\right)&\leq& \frac{\alpha_0}{\delta}\sum_{\frac{1}{|k|^{\tau+1}}< \frac{8\delta\sqrt{\varepsilon_0}}{\alpha}}\frac{1}{|k|^{5\tau+4}|k_{r+1}|}\\
&\leq&\frac{\alpha_0}{\delta}\frac{(8\delta)^4\varepsilon_0^2}{\alpha^4}\sum_{k\in\mathbb{Z}^{r+1}}\frac{1}{|k|^\tau |k_{r+1}|}\\
&\leq& c\varepsilon_0^2\sum_{0\not=k_{r+1}\in\mathbb{Z}}\frac{1}{|k_{r+1}|^{\tau+1}}\\
&\leq& c\varepsilon_0^2.
\end{eqnarray*}
Then
$$\lim_{\varepsilon_0\rightarrow 0}\frac{\mbox{meas}\left(\bigcup_{0\not=k\in\mathbb{Z}^{r+1}}O_k\right)}{\varepsilon_0}=0.$$
Therefore, \eqref{solution1} holds for most sufficiently small $\varepsilon\in(0,\varepsilon_0)$.

Thus, all solutions of  equation \eqref{exp1} are quasi-periodic with frequencies $\Omega=(\omega_1, \cdots,\omega_r, \sqrt{b})$
for most sufficiently small $\varepsilon>0$.\qed


\section*{References}
\bibliographystyle{elsarticle-num}

\end{document}